\documentclass{amsproc}
\usepackage{graphicx}
\usepackage{amssymb}
\usepackage{epstopdf}
\usepackage{amsmath,amsthm,amscd,amssymb}
\usepackage{latexsym}
\usepackage[colorlinks,citecolor=red,pagebackref,hypertexnames=false]{hyperref}
\usepackage{geometry}                
\numberwithin{equation}{section}
\theoremstyle{plain}
\newtheorem{theorem}{Theorem}[section]
\newtheorem{lemma}[theorem]{Lemma}

\theoremstyle{definition}

\newtheorem{case[theorem]}{Case}

\theoremstyle{remark}
\newtheorem{remark}[theorem]{Remark}
\numberwithin{equation}{section}

\begin{document}

\title{\parbox{14cm}{\centering{An analog of the Furstenberg-Katznelson-Weiss theorem on triangles in sets of positive density in finite field geometries}}}

\author{David Covert, Derrick Hart, Alex Iosevich and Ignacio Uriarte-Tuero} 

\address{Department of Mathematics, University of Missouri, Columbia, MO 65211-4100}
\email{covert@math.missouri.edu} 
\email{hart@math.missouri.edu} 
\email{iosevich@math.missouri.edu} 
\email{ignacio@math.missouri.edu}

\thanks{A. Iosevich was supported by the NSF Grant DMS04-56306}

\begin{abstract} We prove that if the cardinality of a subset of the $2$-dimensional vector space over a finite field with $q$ elements is $\ge \rho q^2$, with $\frac{1}{\sqrt{q}}<<\rho \leq 1$, then it contains an isometric copy of $\ge c\rho q^3$ triangles. \end{abstract} 

\maketitle


\section{Introduction}

A classical result due to Furstenberg, Katznelson and Weiss (\cite{FKW90}; see also \cite{B86}) says that if $E \subset {\mathbb R}^2$ has positive upper Lebesgue density, then for any $\delta>0$, the $\delta$-neighborhood of $E$ contains a congruent copy of a sufficiently large dilate of every three point configuration. An example due to Bourgain shows that if the three point configuration in question is an arithmetic progression, then taking a $\delta$-neighborhood is necessary and the result is not otherwise true. However, it seems reasonable to conjecture that if the three point configuration is non-degenerate in the sense that the three points do not lie on the same line, then a set of positive density contains a sufficiently large dilate of this configuration. 

When the size of the point set is smaller than the dimension of ambient Euclidean space, taking a $\delta$-neighborhood is not necessary, as shown by Bourgain in \cite{B86}. He proves that if $E \subset {\Bbb R}^d$ has positive upper density and $\Delta$ is a $k$-simplex with $k<d$, then $E$ contains a rotated and translated image of every large dilate of $\Delta$. The case $k=d$ and $k=d+1$ remain open, however.  See also, for example, \cite{Berg96}, \cite{F81}, \cite{K07}, \cite{TV06} and \cite{Z99} on related problems and their connections with discrete analogs. 

In the geometry of the integer lattice ${\Bbb Z}^d$, related problems have been recently investigated by Akos Magyar in \cite{M06} and \cite{M07}. In particular, he proves in \cite{M07} that if $d>2k+4$ and $E \subset {\Bbb Z}^d$ has positive upper density, then all large (depending on density of $E$) dilates of a $k$-simplex in ${\Bbb Z}^d$ can be embedded in $E$. Once again, serious difficulties arise when the size of the simplex is sufficiently large with respect to the ambient dimension. 

The purpose of this paper is to investigate an analog of this question in finite field geometries. A step in this direction was taken by the second and third listed authors in \cite{HI07}. They prove that if $E \subset {\Bbb F}_q^d$, the $d$-dimensional vector space over the finite field with $q$ elements with $|E| \ge Cq^{d \frac{k-1}{k}+\frac{k-1}{2}}$ and $\Delta$ is a $k$-dimensional simplex, then there exists $\tau \in {\Bbb F}_q^d$ and $O \in SO_d({\Bbb F}_q)$ such that $\tau+O(\Delta) \subset E$. The result is only non-trivial in the range $d \ge \left(^k_2\right)$ as larger simplexes are out of range of the methods used. 

The purpose of this paper is to address the case of triangles in two-dimensional vector spaces over finite fields. Given $E \subset {\mathbb F}_q^2$, define 
$$T_3(E)=\{(x,y,z) \in E \times E \times E\} \ / \sim$$ with the equivalence relation $\sim$ such that $(x,y,z) \sim (x',y',z')$ if there exists $\tau \in {\mathbb F}_q^2$ and 
$O \in SO_2({\mathbb F}_q)$, the set of two-by-two orthogonal matrices over ${\mathbb F}_q$ with determinant $1$, such that 
$$ (x',y',z')=(O(x)+\tau, O(y)+\tau, O(z)+\tau).$$ 

Our main result is the following. 

\begin{theorem} \label{main} Let $E \subset {\mathbb F}_q^2$, and suppose that 
$$ |E| \ge \rho q^2$$ for some $\frac{C}{\sqrt{q}} \leq \rho \leq 1$ with a sufficiently large constant 
$C>0$. Then there exists $c>0$ such that 
$$ |T_3(E)| \ge c\rho q^3.$$ 
\end{theorem} 

In other words, we show that if $E$ has density $ \ge \rho$, then the set of triangles determined by $E$, up to congruence, has density $\ge  c\rho$, where $\rho$ is allowed to depend on $q$ within the parameters given above. 

\begin{remark} Note that in contrast to the Furstenberg-Katznelson-Weiss result (\cite{FKW90}) we do not use dilations. This is natural because there is no order in ${\Bbb F}_q$, so a reasonable analog of proving a result for sufficiently large dilates of a three-point configuration in Euclidean space is proving it for all dilates in finite field geometries. 
\end{remark} 

\begin{remark} Observe that the density condition $|E| \ge \rho q^2$ immediately tells us that the numbers of three-tuples determined by $E$, up to congruence, is 
$$ \ge \frac{\rho^3 q^6}{q \cdot q^2}=\rho^3 q^3,$$ since the size of the translation group is $q^2$ and the size of the rotation group is $q$. Thus our result can be viewed as shaving off two powers of $\rho$ from this (trivial) estimate. It is conceivable that $\rho q^3$ may be replaced by $cq^3$, for some $0<c<1$, or even $(1-o(1))q^3$. 
\end{remark}

\vskip.125in

\subsection{Finite field analog of Bourgain's example} The following variant of Bourgain's Euclidean construction (see \cite{FKW90}) shows that for general subsets of ${\mathbb F}_q^2$ satisfying the density assumption $|E| \ge \rho q^2$ for some $\rho>0$ it is not possible to recover isometric copies of all three point configurations. 

Let $A \subset {\mathbb F}_q^{*}$, the multiplicative group of ${\mathbb F}_q$, such that $|A| \ge \rho q$ for some $\rho>0$ and
$$ 2A+2A-4A \not={\mathbb F}_q.$$ 

We shall give an (easy) construction of such a set at the end of the argument. Let 

$$ E=\bigcup_{t \in A} S_t,$$ where 
$$ S_t=\{x \in {\mathbb F}_q^2: ||x||=t\},$$ with 
$$ ||x||=x_1^2+x_2^2.$$

\vskip.125in 

It is not difficult to check that $|E| \ge \rho q^2$ using the classical fact that a circle in ${\Bbb F}_q^2$ has $q \pm 1$ points. See Lemma \ref{fouriersphere} below. Now consider a three-tuple 

$$ \left(x, y, \frac{x+y}{2} \right)$$ such that 
$$||x-y|| \notin 2A+2A-4A.$$ 

\vskip.125in 

We claim that such a three-tuple cannot be contained in $E$. We shall argue by contradiction. Indeed, the parallelogram law says that 
$$ 2 \left|\left| \frac{x+y}{2} \right|\right|+2\left|\left| \frac{x-y}{2} \right|\right|=||x||+||y||,$$ so 
$$ ||x-y||=2||x||+2||y||-4 \left|\left| \frac{x+y}{2} \right|\right|,$$ which is an element of $2A+2A-4A$. By construction, $||x-y|| \notin 2A+2A-4A$, so we are done. It remains to show that the set $A$ with the desired properties exists. Let $q$ be a large prime number and denote the elements of the corresponding field ${\Bbb F}_q$ by 
$$ \{0,1,2, \dots, q-1\}.$$ 

Let $A$ consists of multiples of $8$ that are less than or equal to $q/32$. This set is clearly of positive density and $2A+2A-4A \not={\mathbb F}_q$ since all of its elements are even as wrap-around is precluded from taking place by the condition that the largest element of $A$ is $\leq q/32$. 

\begin{remark} It is important to note that we do not {\it know} a single example of this type involving a non-degenerate triangle, one whose vertices do not lie on a line. 
\end{remark} 

\vskip.125in 

\section{Proof of the main result (Theorem \ref{main})}

\vskip.125in 

We prove Theorem \ref{main} by reducing it to a statistically more precise statement about hinges. More precisely, we observe that it suffices to show that if $|E| \ge \rho q^2$, then 

\begin{equation} \label{thepoint} |\{(a,b,c) \in {\mathbb F}_q^3: |T_{a,b,c}(E)|>0\}| \ge c\rho q^3, \end{equation} where 

$$ T_{a,b,c}(E)=\{(x,y,z) \in E \times E \times E: ||x-y||=a, ||x-z||=b, ||y-z||=c\},$$ with 

$$ ||x||=x_1^2+x_2^2.$$ 

\vskip.125in 

This follows from the following simple lemma from \cite{HI07}, which we prove at the end of the paper for the sake of completeness. 

\begin{lemma}\label{uptocrap} Let $P$ be a (non-degenerate) simplex with vertices $V_0, V_1, \dots, V_k$, with $V_j \in {\mathbb F}_q^d$. Let $P'$ be another (non-degenerate) simplex with vertices $V'_0, V'_1, \dots, V'_k$. Suppose that
\begin{equation} \label{equalnorm} ||V_i-V_j||=||V'_i-V'_j|| \end{equation} for all $i,j$. Then there exists $\tau \in {\mathbb F}_q^d$ and $O \in SO_d({\mathbb F}_q)$ such that $\tau+O(P)=P'$.
\end{lemma}

The key estimate is the following result about hinges, which is interesting in its own right. 

\begin{theorem} \label{whatwereallyprove} Suppose that $E \subset {\mathbb F}_q^2$ and let $a,b \not=0$. Then 
$$ |\{(x,y,z) \in E \times E \times E: ||x-y||=a, ||x-z||=b\}|={|E|}^3q^{-2}+O(q|E|).$$ 
\end{theorem} 

We can use this result as follows. If $|E|>>q^{\frac{3}{2}}$, then 
$$  |\{(x,y,z) \in E \times E \times E: ||x-y||=a, ||x-z||=b\}|={|E|}^3q^{-2}(1+o(1)).$$ 

By the pigeon-hole principle, there exists $x \in E$ such that 
$$ |\{(y,z) \in E \times E: ||x-y||=a, ||x-z||=b\}| \ge {|E|}^2q^{-2}.$$ 

Suppose that the number of elements of $SO_2({\mathbb F}_q)$ that leave $x$ fixed and keep $(y,z)$ inside the pinned hinge is $ \leq \rho q$. Then recalling our assumption that $$ |E| \ge  \rho q^2 \; ,$$ we get that the number of distinct distances $c$ from $\{y \in E: ||x-y||=a\}$ to $\{z \in E: ||x-z||=b\}$ is at least
$$ {|E|}^2q^{-2}  \; \frac{1}{\rho q} \ge \frac{1}{2} \rho q ,$$ and hence, since there are $(q-1)$ possible choices for $a$ and $b$, \eqref{thepoint} follows. 

If the number of elements of $SO_2({\mathbb F}_q)$ that leave $x$ fixed  and keep $(y,z)$ inside the pinned hinge is $>\rho q$, then both the circle of radius $a$, centered at $x$, and the circle of radius $b$, centered at $x$, contain more than $\rho q$ elements of $E$. The following simple lemma, whose proof is given at the end of this paper, shows that this implies that the number of distinct distances $c$ from 
$$\{y \in E: ||x-y||=a\} \ \text{to} \ \{z \in E: ||x-z||=b\}$$ is $\ge \frac{1}{4} \rho q$, and thus \eqref{thepoint} follows. 

\vskip.125in

\begin{lemma}\label{IntersectingCircles}
Suppose that $a,b,c \not=0$. Let $w=(w_1,w_2) \in \{y \in {\mathbb F}_q^2: ||x-y||=a\}$. 
Consider the set 
$$ I = \{z \in {\mathbb F}_q^2: ||x-z||=b\} \cap \{u \in {\mathbb F}_q^2: ||w-u||=c\} \; .$$

For at least $\frac{q -3}{2}$ different values of $c$, we have that $I \neq \emptyset$.
\end{lemma} 

Thus the proof of Theorem \ref{main} has been reduced to Theorem \ref{whatwereallyprove}. 

\vskip.125in 

\subsection{Fourier analysis used in this paper} Let $\mathbb F _q^d$ be the $d$-dimensional vector space over the finite field ${\mathbb F}_q$. The Fourier transform of a function
$$f: \mathbb F_q^d \rightarrow \mathbb C$$ is given by
$$ \widehat{f}(m)= q^{-d} \sum_{x \in {\mathbb F}_q^d} f(x) \chi(-x \cdot m),$$
where $\chi$ is an additive character on $\mathbb F_q$.

The orthogonality property of the Fourier Transform says that
$$ q^{-d}\sum_{x \in \mathbb F_q^d} \chi(-x \cdot m)=1$$ for $m=(0, \dots, 0)$ and $0$ otherwise. This property yields many of the standard properties of the Fourier Transform. We summarize the basic properties of the Fourier Transform used in this paper as follows.

\begin{lemma} \label{ft} 

Let
$$f,g:\mathbb F_q^d \rightarrow \mathbb C.$$

Then 

$$ \hat f (0, \dots, 0) =q^{-d} \sum_{x \in \mathbb F_q^d} f(x),$$ 
$$ q^{-d} \sum_{x \in \mathbb F_q^d} f(x) \overline{g(x)} =\sum_{m \in \mathbb F_q^d} \hat{f}(m) 
\overline{\hat{g}(m)},$$ and 
$$ f(x) =\sum_{m \in \mathbb F_q^d} \hat{f}(m) \chi(x \cdot m).$$
\end{lemma}

\vskip.25in 

\section{Proof of Theorem \ref{whatwereallyprove}}

\vskip.125in 

Let 

$$ D_a=\{(x,y) \in E \times E: ||x-y||=a\}$$ and identify $D_a$ with its indicator function. We need the following result from \cite{IR07}, proved at the end of this paper for the sake of completeness. 

\begin{theorem} \label{twopoints} Let $E \subset {\mathbb F}_q^d$ and suppose that $t \not=0$. Then 
$$ \sum_{x,y} D_t(x,y)={|E|}^2q^{-1}+O(q^{\frac{d-1}{2}} |E|).$$ 
\end{theorem} 
Let 

$$ S_a=\{x \in {\mathbb F}_q^d: x_1^2+\dots+x_d^2=a\}$$ and identify $S_a$ with its indicator function. 
Now setting $d=2$ and using Lemma \ref{ft} we see that 
$$ |\{(x,y,z) \in E \times E \times E: ||x-y||=a, ||x-z||=b\}|$$
$$=\sum_{x,y,z} D_a(x,y)E(z)S_b(x-z) $$

\begin{equation} \label{fourierexpression}=q^6 \sum_{m} \widehat{D}_a(m,0,0) \widehat{E}(m) \widehat{S}_b(m) . \end{equation} 

We need the following results about the Fourier transform of the sphere which we state in the $d$-dimensional context. 

\begin{lemma} \label{fouriersphere} 
Let  $d \ge 2$ and define $S_b=\{x \in {\mathbb F}_q^d: x_1^2+\dots+x_d^2=b\}.$ Suppose that $b \not=0$ and $m \not=(0,\dots, 0)$. 
\begin{equation} \label{decay} |\widehat{S}_b(m)| \leq 2q^{-\frac{d+1}{2}}. \end{equation}

For any $a \in {\mathbb F}^{*}_q$, 
\begin{equation} \label{size} |S_a|=q^{d-1}+o(q^{d-1}). \end{equation} 
\end{lemma} 

We postpone the proof of the lemma until the end of the paper. In the meantime, we see that Lemma \ref{ft} implies that the expression in (\ref{fourierexpression}) equals 
$$ |D_a| |E||S_b|q^{-2}$$
$$+q^6 \cdot \sum_{m \not=(0,0)} \widehat{D}_a(m, 0,0) \widehat{E}(m) \widehat{S}_b(m)=I+R(a,b).$$ 

In view of Lemma \ref{fouriersphere} and Theorem \ref{twopoints}, 
$$ I={|E|}^3q^{-2}(1+o(1)).$$ 

We have 
$$ R(a,b)=q^6 \sum_{m \not=(0,0)} \widehat{D}_a(m, 0,0) \widehat{E}(m) \widehat{S}_b(m).$$

Using Lemma \ref{fouriersphere} once again and applying Cauchy-Schwartz followed by Lemma \ref{ft}, we see that 

$$ {|R(a,b)|}^2 \leq 4q^{-3} \cdot q^{12} \sum_m {|\widehat{E}(m)|}^2 \cdot \sum_m {|\widehat{D}_a(m,0,0)|}^2$$

\begin{equation} \label{endisnear}=4q^9 \cdot q^{-2} |E| \cdot \sum_m {|\widehat{D}_a(m,0,0)|}^2,
\end{equation} and thus the matter is reduced to bounds for 
$$ \sum_m {|\widehat{D}_a(m,0,0)|}^2.$$ 

By definition, 
$$ \widehat{D}_a(m,0,0)=q^{-4} \sum_{x,y} \chi(-x \cdot m) E(x)E(y)S_a(x-y)$$
$$=q^{-2} \cdot q^{-2} \sum_x \chi(-x \cdot m) E(x)|E \cap S_a(x)|,$$ where 
$$ S_a(x)=\{y \in {\mathbb F}_q^2: ||x-y||=a\}.$$ 

It follows from above and Lemma \ref{ft} that 
$$ \sum_m {|\widehat{D}_a(m,0,0)|}^2=q^{-4} \sum_m {|\widehat{f}(m)|}^2$$
\begin{equation} \label{erection}=q^{-4} \cdot q^{-2} \sum_x {|f(x)|}^2, \end{equation}
where
$$ f(x)=E(x)|E \cap S_a(x)|,$$ and matters are reduced to the estimation of 
$$\sum_x {|f(x)|}^2.$$

\vskip.125in

\begin{lemma} \label{Icantasteit} With the notation above, if $q$ is sufficiently large, then 
$$ \sum_{x \in E} {|E \cap S_a(x)|}^2 \leq 8q|E|.$$ 
\end{lemma} 

\vskip.125in

To prove the result, we write 
$$ |E \cap S_a(x)|=\sum_y E(y) S_a(x-y)$$
$$=|E||S_a|q^{-2}+q^2 \sum_{m \not=(0,0)} \chi(x \cdot m) \widehat{E}(m) \widehat{S}_a(m)$$
$$=A+B(x).$$ 

It is easy to see that plugging in $A$ leads to a better estimate than claimed. Now, 
$$ \sum_{x \in E} {|B(x)|}^2 \leq \sum_{x \in {\mathbb F}_q^2} {|B(x)|}^2$$
$$ =q^4 \sum_{m,m' \not=(0,0)} \widehat{E}(m) \overline{\widehat{E}(m')} \widehat{S}_a(m) \overline{\widehat{S}_a(m')} \sum_x \chi(x \cdot (m-m'))$$
$$=q^6 \sum_{m \not=(0,0)} {|\widehat{E}(m)|}^2 {|\widehat{S}_a(m)|}^2,$$ and by Lemma \ref{fouriersphere} this quantity is 
$$ \leq 4q^6 \cdot q^{-3} \sum_{m \not=(0,0)} {|\widehat{E}(m)|}^2,$$ which, by Lemma \ref{ft} is 
$$ \leq 4q|E|,$$ as claimed. 

Plugging everything back into (\ref{erection}) and then (\ref{endisnear}), we see that 
\begin{equation} \label{orgasm} 
{|R(a,b)|}^2 \leq 32 q^9 \cdot q^{-2} \cdot |E| \cdot q^{-6} \cdot q|E| \end{equation}
$$=32q^2 {|E|}^2.$$ 

Recall that 
$$ |\{(x,y,z) \in E \times E \times E: ||x-y||=a, ||x-z||=b\}|=I+R(a,b),$$ where 
$$ I={|E|}^3q^{-2}(1+o(1))$$ and 
$$ |R(a,b)| \leq 8 q|E|$$ for $q$ sufficiently large by (\ref{orgasm}) above. 

It follows that 
$$  |\{(x,y,z) \in E \times E \times E: ||x-y||=a, ||x-z||=b\}|={|E|}^3q^{-2}(1+o(1))$$ if $|E|>Cq^{\frac{3}{2}}$ with a sufficiently large constant $C>0$, as claimed. 

\vskip.125in 

\section{Proof of Theorem \ref{twopoints}}

\vskip.125in 

We have 
$$ \sum_{x,y} D_t(x,y)=\sum_{x,y}E(x)E(y)S_t(x-y).$$ 
Applying Fourier inversion to the sphere,
\begin{align*}
\sum_{x,y} D_t(x,y)=\sum_{x,y} E(x)E(y)\sum_m \widehat{S}_{t}(m)\chi(m\cdot(x-y))
\\
&=q^{2d} \sum_m {|\widehat{E}(m)|}^2 \widehat{S}_{t}(m)
\\
&={|E|}^2 \cdot q^{-d} \cdot |S_{t}|+q^{2d} \sum_{m \not=(0,
\dots, 0)} {|\widehat{E}(m)|}^2 \widehat{S}_{t}(m)
\\
&=M+R.\,
\end{align*}

By Lemma \ref{fouriersphere},
$$ M =\frac{{|E|}^2}{q}+|E|^2o(q^{-1}),$$ and using Lemma \ref{fouriersphere} once
again, along with Lemma \ref{ft}, we have 
\begin{align*}
|R| &\leq 2 q^{2d} \cdot q^{-\frac{d+1}{2}} \cdot \sum_m {|\widehat
{E}(m)|}^2
\\
&=2q^{\frac{d-1}{2}} \cdot |E|.
\end{align*}

This completes the proof.

\section{Proof of Lemma \ref{fouriersphere}}

\vskip.125in

For any $l\in{\mathbb F}^d_q$, we have
\begin{equation} \label{sphereparade}
\begin{array}{llllll} \widehat{S}_t(l)&=&
q^{-d} \sum_{x \in {\mathbb F}^d_q} q^{-1} \sum_{j \in {\mathbb F}_q} \chi(
j(\|x\|-t)) \chi( -  x \cdot l)\\ \hfill \\&=&q^{-1}\delta(l) +
q^{-d-1} \sum_{j \in {\mathbb F}^{*}_q} \chi(-jt) \sum_{x}
\chi( j\|x\|) \chi(- x \cdot l),\\
\end{array}\end{equation}where the notation $\delta(l)=1$ if $l=(0\ldots,0)$ and $\delta(l)=0$ otherwise. 

Now
$$\widehat{S}_t(l)=q^{-1}\delta(l)+   Q^d q^{-\frac{d+2}{2}} \sum_{j \in {\mathbb F}^{*}_q}
\chi\left(\frac{\|l\|}{4j}+jt\right)\eta^d(-j).$$
In the last line we have completed the square, changed $j$ to
$-j$, and used $d$ times the Gauss sum equality (see e.g. \cite{LN97}) 
\begin{equation} \sum_{c\in {\mathbb F}_q} \chi(jc^2) = \eta(j)\sum_{c\in{\mathbb F}_q}\eta(c)\chi(c)=\eta(j)\sum_{c\in{\mathbb F}_q^*}\eta(c)\chi(c) =Q\sqrt{q}\,\eta(j), \label{gauss}\end{equation} where the constant $Q$ equals $\pm1$ or $\pm i$, depending on $q$, and $\eta$ is the quadratic multiplicative character (or the Legendre symbol) of ${\mathbb F}_q^*$.

The conclusion to the first and second parts of Lemma \ref{fouriersphere} now follows from standard Gauss sum estimates (see e.g. \cite{LN97}) and the following classical estimate due to A. Weil (\cite{W48}).

\begin{theorem} \label{kloosterman} Let
$$ K(a)=\sum_{s \not=0} \chi(as+s^{-1}) \psi(s), $$ where $\psi$ is a multiplicative character on ${\mathbb F}_q^{*}$. Then, if $a \not=0$,
$$ |K(a)| \leq 2 \sqrt{q}.$$ 
\end{theorem}

\vskip.125in 

\section{Proof of Lemma \ref{uptocrap}}

\vskip.125in

Let $\pi_r(x)$ denote the $r$th coordinate of $x$. Taking translations into account, we may assume that $V_0=(0, \dots, 0)$. We may also assume that $V_1, \dots, V_k$ are contained in ${\mathbb F}_q^k$. The condition (\ref{equalnorm}) implies that
\begin{equation} \label{dotproduct} \sum_{r=1}^k \pi_r(V_i) \pi_r(V_j)=
\sum_{r=1}^k \pi_r(V_i') \pi_r(V_j'). \end{equation}

Let $T$ be the linear transformation uniquely determined by the condition
$$ T(V_i)=V'_i.$$

In order to prove that $T$ is orthogonal, it suffices to show that
$$ ||Tx||=||x||$$ for any $x \not=(0, \dots, 0)$. We give this (standard) reduction below for the sake of reader's convenience. 

Since $V_j$s form a basis, by assumption, we have
$$ x=\sum_i t_i V_i, $$ so it suffices to show that
$$ ||x||=\sum_r \sum_{i,j} t_i t_j \pi_r(V_i) \pi_r(V_j)$$
$$=\sum_r \sum_{i,j} t_i t_j \pi_r(V'_i) \pi_r(V'_j)=||Tx||,$$ which follows
immediately from (\ref{dotproduct}).

Observe that we used the fact that orthogonality of $T$, the condition that
$T^t \cdot T=I$ is equivalent to the condition that $||Tx||=||x||$. To see
this observe that to show that $T^t \cdot T=I$ it suffices to show that
$T^tTx=x$ for all non-zero $x$. This, in turn, is equivalent to the
statement that
$$ <T^tTx,x>=||x||,$$ where
$$ <x,y>=\sum_{i=1}^k x_iy_i.$$

Now,
$$ <T^tTx,x>=<Tx, Tx>$$ by definition of the transpose, so the stated
equivalence is established. This completes the proof of Lemma
\ref{uptocrap}.

\vskip.125in

\section{Proof of Lemma \ref{IntersectingCircles}}

\vskip.125in

After a translation, we may assume without loss of generality that $x = (0,0)$. We are looking for solutions $(s,t)$ to the system of equations
\begin{eqnarray}\label{e.IntersectingCircles1}
\left\{ \begin{array}{l}
(s-w_1)^2 + (t-w_2)^2 =c \\ s^2 + t^2 =b \end{array}  \right. \end{eqnarray} with the assumption that $w_1^2 + w_2^2 = a$. Then \eqref{e.IntersectingCircles1} is equivalent to
\begin{eqnarray}\label{e.IntersectingCircles2} \left\{ \begin{array}{l} w_1 \cdot s + w_2 \cdot t =\frac{a+b-c}{2} \\ s^2 + t^2 =b \end{array}   \right. \end{eqnarray}

Now $w_1$ and $w_2$ cannot be simultaneously zero since $a \neq 0$. If $w_1 \neq 0$, from the first equation in \eqref{e.IntersectingCircles2} we get that
\begin{equation}\label{SolvingForS}
s = \frac{1}{w_1} \left\{ \frac{a+b-c}{2} - w_2 \cdot t \right\} \; ,
\end{equation} which substituted into the second equation in \eqref{e.IntersectingCircles2} gives
\begin{equation}\label{e.IntersectingCircles3} \left\{ \left( \frac{w_2}{w_1} \right)^2  +1 \right\} t^2 - \left\{ \frac{a+b-c}{w_1} \; \frac{w_2}{w_1} \right\} t + \left( \frac{a+b-c}{2w_1}  \right)^2 -b =0 \; . \end{equation}

(If $w_1 = 0$ so that $w_2 \neq 0$, the resulting equation is the same as \eqref{e.IntersectingCircles3} but interchanging the roles of $w_1$ and $w_2$ among themselves and the roles of $s$ and $t$ among themselves.)

However, notice now that the condition $w_2 = \pm i \; w_1$ is incompatible with the hypothesis that $w_1^2 + w_2^2 = a \neq 0$. Consequently, $\left( \frac{w_2}{w_1} \right)^2 \not= -1 $, and hence the equation \eqref{e.IntersectingCircles3} has at most 2 solutions. We still have to prove that  the equation \eqref{e.IntersectingCircles3} has indeed a solution under our hypotheses. The discriminant of equation \eqref{e.IntersectingCircles3} is

\begin{eqnarray}\label{Discriminant.e.IntersectingCircles3} \Delta  =  \left\{ \frac{a+b-c}{w_1} \; \frac{w_2}{w_1} \right\}^2 -4\left\{ \left( \frac{w_2}{w_1} \right)^2  +1 \right\} \left[ \left( \frac{a+b-c}{2w_1}  \right)^2 -b \right] & = & \frac{ 4b\left[w_1^2 + w_2^2 \right] -\left( a+b-c \right)^2}{w_1^2}  =  \nonumber \\ & = & \frac{4ab -\left( a+b-c \right)^2}{w_1^2}\; . \end{eqnarray}

Hence equation \eqref{e.IntersectingCircles3} has a solution precisely when there exists a $k \in \mathbb F_q$ such that 
\begin{equation}\label{ConditionOnSquareRootDiscriminant.e.IntersectingCircles3} 4ab -\left( a+b-c \right)^2 = k^2 \; , \end{equation} which happens precisely when $c$ is of the form
\begin{equation}\label{SolvingForCInConditionOnSquareRootDiscriminant.e.IntersectingCircles3} c = a+b \pm \sqrt{4ab-k^2} \; , \end{equation} i.e. whenever there exists a $\tau \in \mathbb F_q$ such that
\begin{equation}\label{ConditionOfSumOfSquaresAfterSolvingForCInConditionOnSquareRootDiscriminant.e.IntersectingCircles3} k^2 + \tau^2 = 4ab \; .
\end{equation}

We now repeat, for the convenience of the reader, the well-known known argument that every element of $\mathbb F_q$ is a sum of $2$ squares. Namely, by \eqref{gauss}, and recalling that $\eta^2(t) = 1$ and that $\sum_{t \in \mathbb F_q} \chi(-4abt) =0$ if $t \neq 0$, and that $Q = \pm 1$ or $\pm i$ depending on $q$,

\begin{eqnarray}\label{CountingSumsOfSquares1} |\{ (k, \tau): k^2 + \tau^2 = 4ab  \}| & = & \frac{1}{q} \sum_{t \in \mathbb F_q} \; \sum_{k, \tau \in \mathbb F_q} \chi(t(k^2 + \tau^2 - 4ab)) = \nonumber \\
& = & q + \frac{1}{q} \sum_{t \neq 0} \; \sum_{k, \tau \in \mathbb F_q} \chi(t(k^2 + \tau^2 - 4ab)) = \nonumber \\ & = & q + \frac{1}{q} \sum_{t \neq 0} \chi(-4abt) \; Q^2 \; q \,\eta^2(t) = q -Q^2 \; .
\end{eqnarray}

Hence equation \eqref{e.IntersectingCircles3} has a solution for at least $\frac{q \pm 1}{2}$ different values of $c$, since by \eqref{SolvingForCInConditionOnSquareRootDiscriminant.e.IntersectingCircles3} and 
\eqref{ConditionOfSumOfSquaresAfterSolvingForCInConditionOnSquareRootDiscriminant.e.IntersectingCircles3} each value of $c$ for which equation \eqref{e.IntersectingCircles3} has a solution corresponds precisely to one value of $\tau$, and each value of $\tau$ is accounted for at most twice in \eqref{CountingSumsOfSquares1} since for each such value of $\tau$, there are at most $2$ values of $k$ satisfying \eqref{ConditionOfSumOfSquaresAfterSolvingForCInConditionOnSquareRootDiscriminant.e.IntersectingCircles3}.

Since it is conceivable (depending on the value of $q$) that $c=0$ yields a solution to \eqref{e.IntersectingCircles1}, accounting for that possibility, we can assert that \eqref{e.IntersectingCircles1} has a solution for at least $\frac{q -3}{2}$ different values of $c \neq 0$.


\vskip.125in

\newpage

\begin{thebibliography}{4}

\bibitem{B86} J. Bourgain, {\it A SzemerŽdi type theorem for sets of positive density}, Israel J. Math. \textbf{54} (1986), no. 3, 307-331.

\bibitem{Berg96} V. Bergelson, {\it Ergodic Ramsey theory Ñ an update}, Ergodic Theory of Zd -Actions (Warwick, 1993 Ð 1994) (M. Pollicott and K. Schmidt, eds.), London Math. Soc. Lecture Note Ser., \textbf{228}, Cambridge Univ. Press, Cambridge, (1996). 

\bibitem{F81} H. Furstenberg, {\it Recurrence in ergodic theory and combinatorial number theory}, M. B. Porter Lectures, Princeton Univ. Press, Princeton, NJ, (1981). 

\bibitem{FKW90} H. Furstenberg, Y. Katznelson, and B. Weiss, {\it Ergodic theory and configurations in sets of positive density} Mathematics of Ramsey theory, 184-198, Algorithms Combin., 5, Springer, Berlin, (1990).

\bibitem{HI07} D. Hart and A. Iosevich, {\it Ubiquity of simplices in subsets of vector spaces over finite
fields}, Analysis Mathematika, \textbf{34}, (2007). 

\bibitem{IR07} A. Iosevich and M. Rudnev, {\it Erd\H os distance problem in vector spaces over finite fields}. Trans. Amer. Math. Soc. (2007).

\bibitem{K07} B. Kra, {\it Ergodic methods in additive combinatorics}, Centre de Recherches Mathematiques Proceedings and Lecture Notes (2007). 

\bibitem{LN97} R. Lidl and H. Niederreiter, {\it Finite fields}, Cambridge University Press, (1997). 

\bibitem{M06} A. Magyar, {\it On distance sets of large sets of integers points}, Israel Math J. (to appear) (2006). 

\bibitem{M07} A. Magyar, {\it $k$-point configurations in sets of positive density of ${\Bbb Z}^n$}, Duke Math J. (to appear), (2007). 

\bibitem{TV06} T. Tao and V. Vu. {\it Additive Combinatorics}. Cambridge
University Press, 2006.

\bibitem{W48} A. Weil, {\it On some exponential sums}, Proc. Nat. Acad. Sci. U.S.A. \textbf{34} (1948) 204-207.

\bibitem{Z99} T. Ziegler, {\it An application of ergodic theory to a problem in geometric Ramsey theory}, Israel Journal of Math. \textbf{114} (1999) 271-288.

\end{thebibliography}
\end{document}